\documentclass[a4paper,10pt]{article}
\usepackage{mathtext}
\usepackage[T1,T2A]{fontenc}
\usepackage[cp1251]{inputenc}
\usepackage[english]{babel}
\usepackage{amsmath}
\usepackage{amsfonts}
\usepackage{amssymb}
\usepackage{mathrsfs}
\usepackage{amsthm}
\usepackage{titling}
\usepackage{multicol}
\usepackage{graphicx}
\usepackage{wrapfig}
\usepackage{appendix}

\usepackage{color}
\usepackage{euscript}

\newcommand{\Bell}{\boldsymbol{B}}

\newcommand{\Class}{\boldsymbol{A}}

\newcommand{\BMO}{\mathrm{BMO}}
\newcommand{\eps}{\varepsilon}
\DeclareMathOperator{\cl}{cl}

\newcommand{\av}[2]{\langle {#1}\rangle_{{}_{#2}}}

\renewcommand{\leq}{\leqslant}
\renewcommand{\geq}{\geqslant}

\newcommand{\BG}{\mathfrak{B}}

\newtheorem{Le}{Lemma}[section]
\newtheorem{Def}[Le]{Definition}

\newtheorem{Th}[Le]{Theorem}

\newtheorem{Conj}[Le]{Conjecture}

\newtheorem{Prob}[Le]{Problem}
\numberwithin{equation}{section}

\begin{document}
\author{P. Ivanisvili\thanks{The author visited Hausdorff Research Institute for Mathematics (HIM) in the framework of the Trimester Program ``Harmonic Analysis and Partial Differential Equations''. He thanks HIM for the hospitality.} \and N. N. Osipov\thanks{This research is supported by the Chebyshev Laboratory  (Department of Mathematics and Mechanics, St. Petersburg State University)  under RF Government grant 11.G34.31.0026.}\;\thanks{Supported by RFBR grant no. 14-01-31163.}\;\thanks{Supported by RFBR grant no. 14-01-00198.}\;\thanks{Supported by ERCIM ``Alain Bensoussan'' Fellowship.} \and D.~M.~Stolyarov\thanksmark{2}\;\thanks{Research was supported by JSC "Gazprom Neft".}\;\thanksmark{4} \and V. I. Vasyunin\thanks{Supported by RFBR grant no. 14-01-00748a.} \and P.~B.~Zatitskiy\thanksmark{2}\;\thanksmark{6}\;\thanks{Supported by the SPbSU grant no. 6.38.223.2014.}\;\;\,\thanks{Supported by the Dinasty Foundation.}}
\title{Sharp estimates of integral functionals on classes of functions with small mean oscillation}

\maketitle
\begin{abstract}
We unify several Bellman function problems treated in~\cite{BR,DW,IOSVZ,IOSVZShortReport,ISVZ,LSSVZ,Osekowski2,Reznikov,Slavin,SSV,SV,SV2,Vasyunin,Vasyunin2,Vasyunin3,Vasyunin4,Vasyunin5,VV,VV1} into one setting. For that purpose we define a class of functions that have, in a sense, small mean oscillation 
(this class depends on two convex sets in~$\mathbb{R}^2$).
We show how the unit ball in the~$\BMO$ space, or a Muckenhoupt class, or a Gehring class can be described in such a fashion. Finally, we consider a Bellman function problem on these classes, discuss its solution and related questions.   
\end{abstract}

Since Slavin~\cite{Slavin} and Vasyunin~\cite{Vasyunin} proved the sharp form of the John--Nirenberg inequality (see~\cite{SV}), there have been many papers where similar principles are used to prove sharp estimates of this kind. However, there is no theory or even a unifying approach; moreover, the class of problems to which the method can be applied has not been described yet. There is a portion of heuristics in the folklore that is each time applied to a new problem in a very similar manner. The first attempt to build a theory (at least for $\BMO$) was made in~\cite{SV2}, then the theory was developed in the paper~\cite{IOSVZ} (see the short report~\cite{IOSVZShortReport} also). We would also like to draw the reader's attention to the forthcoming paper~\cite{ISVZ}, which can be considered as a description of the theory  for the~$\BMO$ space in a sufficient generality. Problems of this kind were considered not only in~$\BMO$, but in Muckenhoupt classes, Gehring classes, etc 
(see~\cite{BR, DW, Reznikov, Slavin2, Vasyunin2, Vasyunin3}). In this short note, we define a class of functions and an extremal problem on it that includes all the problems discussed above. We believe that the unification we offer gives a strong basis for a theory that will distinguish a certain class of problems to which the method is applicable in a direct way. In Section~\ref{S1} we state the problem and discuss related questions. Section~\ref{S2} contains a detailed explanation how our classes of functions include the unit ball in~$\BMO$ as well as the ``unit balls'' in Muckenhoupt classes and Gehring classes. Finally, in Section~\ref{S3} we give hints to the solution of the problem (as the reader may expect looking at previous papers, it is rather lengthy and technical, so we omit a description of the solution, but concentrate on an analogy with the case of~$\BMO$ considered in~\cite{IOSVZ, ISVZ, Vasyunin5}).

\section{Setting}\label{S1}
Let~$\Omega_0$ be a non-empty open strictly convex subset of~$\mathbb{R}^2$ and let~$\Omega_1$ be open strictly convex subset of~$\Omega_0$. We define the domain~$\Omega$ 
as~$\cl( \Omega_0 \setminus  \Omega_1)$ (the word ``domain'' comes from ``domain of a function''; the symbol~$\cl$ denotes the closure) and the class~$\Class_{\Omega}$ of summable~$\mathbb{R}^2$-valued functions on an interval~$I \subset \mathbb{R}$ as follows:
\begin{equation}\label{AnalyticClass}
\Class_{\Omega} = 
\big\{\varphi \in L^1(I,\mathbb{R}^2) \;\mid\; \varphi(I) \subset \partial \Omega_0 \quad\mbox{and}\quad \forall \;\mbox{subinterval}\; J\subset I \quad \av{\varphi}{J} \notin \Omega_1\big\}.
\end{equation}
Here~$\av{\varphi}{J} = \frac{1}{|J|}\int_J \varphi(s)\,ds$ is the average of~$\varphi$ over~$J$.
In Section~\ref{S2} we show how the unit ball in~$\BMO$ as well as the ``unit balls'' in Muckenhoupt and Gehring classes can be represented in the form~\eqref{AnalyticClass}. Let~$f\colon \partial \Omega_0 \to \mathbb{R}$ be a bounded from below Borel measurable locally bounded function. We are interested in sharp bounds for the expressions of the form
$\av{f(\varphi)}{I}$, where $\varphi \in \Class_{\Omega}$.

Again, in Section~\ref{S2} we explain how the John--Nirenberg inequality or other inequalities of harmonic analysis can be rewritten as estimations of such an expression. The said estimates are delivered by the corresponding Bellman function
\begin{equation}\label{BellmanFunction}
\Bell_{\Omega,f}(x) = \sup\big\{\av{f(\varphi)}{I}\,\big|\,\, \av{\varphi}{I} = x,\,\,\varphi \in \Class_{\Omega}\big\} .
\end{equation}
\begin{Prob}\label{MainProblem}
Given a domain~$\Omega$ and a function~$f$\textup, calculate the function~$\Bell_{\Omega,f}$.
\end{Prob}
As it has been said in the abstract, the particular cases of this problem were treated in 
the papers~\cite{BR,DW,IOSVZ,IOSVZShortReport,ISVZ,LSSVZ,Osekowski2,Reznikov,Slavin,SSV,SV,SV2,Vasyunin,Vasyunin2,Vasyunin3,Vasyunin4,Vasyunin5,VV,VV1} (see Section~\ref{S2} for a detailed explanation). The main reason for Problem~\ref{MainProblem} to be solvable (and it has been heavily used in all the preceeding work) is that the function~$\Bell$ enjoys good properties.
\begin{Def}\label{LocalConcavity}
Let~$\omega$ be a subset of~$\mathbb{R}^d$. We call a function~$G\colon w \to \mathbb{R}\cup\{+\infty\}$ locally concave on~$\omega$ if for every segment~$\ell \subset \omega$ the restriction~$G\big|_{\ell}$ is concave.
\end{Def}
Define the class of functions on~$\Omega$:
\begin{equation}\label{Lambdaclass}
\Lambda_{\Omega,f} = \Big\{G \colon \Omega\to \mathbb{R}\cup\{+\infty\}
\,\Big|\;
G\hbox{ is locally concave on~$\Omega$,} \quad \forall x \in \partial \Omega_0 \quad G(x) \geq f(x)\Big\}.
\end{equation}
The function~$\BG_{\Omega,f}$ is given as follows:
\begin{equation}\label{MinimalLocallyConcave}
\BG_{\Omega,f}(x) = \inf_{G \in \Lambda_{\Omega,f}} G(x).
\end{equation}
\begin{Conj}\label{DualityConjecture}
$\Bell_{\Omega,f} = \BG_{\Omega,f}$.
\end{Conj}
In particular, the conjecture states that the Bellman function is locally concave (because the function~$\BG_{\Omega,f}$ is). 
\begin{Prob}\label{DualityProblem}
Prove Conjecture~\textup{\ref{DualityConjecture}} in adequate generality.
\end{Prob}
Though it may seem that one should solve Problem~\ref{DualityProblem} before turning to Problem~\ref{MainProblem}, it is not really the case. All the preceeding papers used Conjecture~\ref{DualityConjecture} as an assumption that allowed the authors to guess~$\Bell$, then prove that this function was the Bellman function indeed, and only then verify Conjecture~\ref{DualityConjecture} for~$\Omega$ and~$f$ chosen. However, to treat Problem~\ref{DualityProblem} in itself, one has to invent a different approach, see Section~\ref{S3}.

We note that one should impose some additional conditions on~$\Omega$ and~$f$ to provide a solution to the problems. We postpone the detailed discussion of this to Section~\ref{S3} and pass to examples.

\section{Examples}\label{S2}
From now on, we follow the agreement: if~$g\colon \mathbb{R} \to \mathbb{R}^2$ is some fixed parametrization of~$\partial\Omega_0$, then the function~$f(g)\colon \mathbb{R} \to \mathbb{R}$ is denoted by~$\tilde{f}$.

\paragraph{The~$\BMO$ space.} We consider the~$\BMO$ space with the quadratic seminorm. Let~$\eps$ be a positive number. Set~$\Omega_0 = \{x\in\mathbb{R}^2\mid x_1^2 < x_2\}$ and~$\Omega_{1} = \{x\in\mathbb{R}^2\mid  x_1^2 + \eps^2 < x_2\}$. A function $$\varphi = (\varphi_1,\varphi_2)\colon I \to \partial\Omega_0$$ belongs to the class~$\Class_{\Omega}$ if and only if its first 
coordinate $\varphi_1$ belongs to~$\BMO_{\eps}$ (the ball of radius~$\eps$ in~$\BMO$). Indeed, for any~$t \in I$ we have~$\varphi_2(t) = \varphi_1^2(t)$, therefore, the condition~$\av{\varphi}{J}\notin\Omega_1$ can be rewritten as
\begin{equation*}
\av{\varphi_1^2}{J} \leq \av{\varphi_1}{J}^2 + \eps^2,
\end{equation*}
which is the same as
\begin{equation}\label{AlmostBMO}
\av{\big(\varphi_1 - \av{\varphi_1}{J}\big)^2}{J} \leq \eps^2.
\end{equation}
Now we see that the class~$\Class_{\Omega}$ corresponds to~$\BMO_{\eps}$. The Bellman function~\eqref{BellmanFunction} estimates the functional~$\av{\tilde{f}(\varphi_1)}{I}$. The solution of Problem~\ref{MainProblem} with~$\tilde{f}(t) = e^{\lambda t}$ leads to the John--Nirenberg inequality in its integral form, the case~$\tilde{f}(t) = \chi_{(-\infty,-\lambda]\cup[\lambda,\infty)}(t)$ corresponds to the weak form of the John--Nirenberg inequality, and the case~$f(t) = |t|^p$ leads to equivalent defintions of~$\BMO$. We address the reader to the paper~\cite{IOSVZ} for a detailed discussion. This case is the subject of study for the papers~\cite{IOSVZ, ISVZ, LSSVZ, Osekowski2, SV, SV2, Vasyunin4, Vasyunin5}.

\paragraph{Classes~$A_{p_1,p_2}$.} Let~$p_1$ and~$p_2$,~$p_1 > p_2$, be real numbers and let~$Q \geq 1$. Suppose
\begin{equation*}
\Omega_0 = \{x \in\mathbb{R}^2\mid  x_1,x_2 > 0,\,\, x_2^{\frac{1}{p_2}} < x_1^{\frac{1}{p_1}}\} \quad \hbox{and}\quad \Omega_1 = \{x \in\mathbb{R}^2\mid  x_1,x_2 > 0,\,\, Q x_2^{\frac{1}{p_2}} < x_1^{\frac{1}{p_1}}\}. 
\end{equation*}
If a function~$\varphi$ belongs to the class~$\Class_{\Omega}$, then its first coordinate~$\varphi_1$ belongs to the so-called~$A_{p_1,p_2}$ class. The ``norm'' in this class is defined as
\begin{equation}\label{Apnorm}
[\psi]_{A_{p_1,p_2}} = \sup\limits_{J \subset I} \,\,\,\av{\psi^{p_1}}{J}^{\frac{1}{p_1}}\av{\psi^{p_2}}{J}^{-\frac{1}{p_2}},
\end{equation}
where the supremum is taken over all subintervals of~$I$.
These classes were introduced in~\cite{Vasyunin3}. If~$p \in (1,\infty)$, then~$A_{1,-\frac{1}{p-1}} = A_p$, where~$A_p$ stands for the classical Muckenhoupt class. The limiting cases~$A_1$ and~$A_{\infty}$ also fit into this definition (with Hruschev's ``norm'' on~$A_{\infty}$). When~$p_2 = 1$ and~$p_1 > 1$, the class~$A_{p_1,p_2}$ coincides with the so-called Gehring class (see~\cite{Korenovskii2} or~\cite{KS}). One can see that the functions in the Gehring class are exactly those that satisfy the reverse H\"older inequality. Sometimes, the Gehring class is called the reverse-H\"older class. Estimates of integral functionals as provided by the Bellman function~\eqref{BellmanFunction} lead to various sharp forms of the reverse H\"older inequality, see~\cite{Vasyunin3}. These cases were treated in the papers~\cite{BR, DW, Reznikov, Vasyunin2, Vasyunin3}.

\paragraph{Reverse Jensen classes.} These classes were introduced in~\cite{Korenovskii2}. Let~$\Phi\colon \mathbb{R}_+ \to \mathbb{R}_+$ be a convex function. Let~$Q > 1$. Consider the class of functions~$\psi\colon I \to \mathbb{R}_+$ such that
\begin{equation*}
\forall J \subset I \quad \av{\Phi(\psi)}{J} \leq Q \Phi(\av{\psi}{J}).
\end{equation*}
Surely, both a Muckenhoupt class and a Gehring class can be described as certain Reverse Jensen classes. The corresponding domain is~$\{x \in \mathbb{R}^2 \mid x_1,x_2 \geq 0,\,\, \Phi(x_1) \leq x_2 \leq Q \Phi(x_1)\}$. Consult a very recent paper~\cite{Slavin2}, where the Bellman function on the domain~$\{x\in \mathbb{R}^2\mid e^{x_1} \leq x_2 \leq Ce^{x_1}\}$,~$C > 1$, provides sharp constants in the John--Nirenberg inequality for the~$\BMO$ space equipped with the~$L^p$-type seminorm.

\section{Hints to solutions}\label{S3}
First, we note that strict convexity of~$\Omega_0$ implies the fact that~$\Bell(x) = f(x)$ for~$x\in \partial \Omega_0$.
Second, we need~$\Omega$ to fulfill several assumptions that all the domains listed in Section~\ref{S2} do satisfy.
\begin{align}
\label{FirstCondition} &1.\hbox{ The domains~$\Omega_0$ and~$\Omega_1$ are unbounded.}\\
\label{SecondCondition} &2.\hbox{ The boundary of~$\Omega_1$ is~$C^2$-smooth.}\\
\label{ThirdCondition} &3.\hbox{ Every ray inside~$\Omega_0$ can be translated to belong to~$\Omega_1$ entirely.}
\end{align}
The first two conditions are technical in a sense, the third one is essential, since (under assumption~\eqref{FirstCondition}) it is equivalent to the fact that for any~$x \in \Omega$ there exists a function~$\varphi \in \Class_{\Omega}$ such that~$\av{\varphi}{I} = x$ (i.e. the supremum in formula~\eqref{BellmanFunction} is taken over a non-empty set). Now we are ready to present a solution of Problem~\ref{DualityProblem}.
\begin{Th}\label{DualityTheorem}
Let the domain~$\Omega$ satisfy the 
conditions~\textup{\eqref{FirstCondition}, \eqref{SecondCondition}, \eqref{ThirdCondition}}. 
If the function~$f$ is bounded from below\textup, then~$\BG_{\Omega,f} = \Bell_{\Omega,f}$.
\end{Th}
The condition that~$f$ is bounded from below is not necessary. However, we note that without this condition the extremal problem in formula~\eqref{BellmanFunction} is not well posed (the integral of~$f(\varphi)$ may be not well defined). In~\cite{SZ} the reader can find the proof of Theorem~\ref{DualityTheorem} for the case~$\cl\Omega_1 \subset \Omega_0$ as well as its analog where~$f$ can be unbounded from below.

To solve Problem~\ref{MainProblem}, we need to consider even more restrictive conditions, we introduce some notation for that purpose. Choose~$g=(g_1,g_2)\colon \mathbb{R} \to \mathbb{R}^2$ to be a continuous parametrization of~$\partial\Omega_0$; let the domain~$\Omega$ lie on the left of this oriented curve. For any number~$u \in \mathbb{R}$ we draw two tangents from the point~$g(u)$ to the set~$\Omega_1$; by a tangent we mean not a line, but a segment connecting~$g(u)$ with the tangency point. We denote the lenghts of the left and the right tangents by~$\ell_{\mathrm{L}}(u)$ and~$\ell_{\mathrm{R}}(u)$ correspondingly (the left tangent lies between the right one and~$g'$, see~\cite{IOSVZ} for explanations about this notation). 
\begin{align}
\label{OneMoreFirstCondition} &1. \hbox{ The boundaries~$\partial\Omega_0$ and~$\partial\Omega_1$ are~$C^3$-smooth curves, the function~$f$ is~$C^3$-smooth.}\\
\label{OneMoreSecondCondition} &2. \hbox{ The curve~$\gamma(t)=\big(g_1(t), g_2(t),\tilde{f}(t)\big) \subset \mathbb{R}^3$ changes the sign of its torsion only a finite number of times.}\\ 
\label{OneMoreThirdCondition} &3. \hbox{ The integrals~$\int_{-\infty}^0 \frac{1}{\ell_{\mathrm{R}}}$ and~$\int_{0}^{+\infty} \frac{1}{\ell_{\mathrm{L}}}$ diverge.}
\end{align}
In Condition~\eqref{OneMoreThirdCondition} the integration is with respect to the natural parametrization of the curve~$\partial\Omega_1$, where the functions $\ell_{\mathrm{R}}$ and $\ell_{\mathrm{L}}$ are considered as the functions of their tangency points lying on $\partial\Omega_1$. For the case where~$g(t) = (t,t^2)$ treated in~\cite{IOSVZ}, Condition~\eqref{OneMoreSecondCondition} turns into ``the function~$\tilde{f}'''$ changes its sign only a finite number of times''; this is exactly the regularity condition we used in~\cite{ISVZ}. The last Condition~\eqref{OneMoreThirdCondition} is more mysterious, we believe that our considerations may work without it.

We also need a summability assumption for the function~$f$. Let~$\alpha_{\mathrm{R}}(u)$ denote the oriented angle between the right tangent at the point~$u$ and the vector~$(1,0)$, let~$\alpha_{\mathrm{L}}(u)$ denote the oriented angle between the left tangent at the point~$u$ and the vector~$(1,0)$. Then, the summability condition requires the bulky integral
\begin{equation}\label{BulkyIntegral}
\int\limits_{-\infty}^t \exp{\bigg(\int\limits_{\tau}^t \frac{g_1'}{\ell_{\mathrm{R}}\cos(\alpha_{\mathrm{R}})}\bigg)}\, \frac{\tan(\alpha_{\mathrm{R}}(\tau)) g_1'(\tau)-g_2'(\tau)}{(g_1'(\tau)g_2''(\tau)-g_2'(\tau)g_1''(\tau))^2}\,
\begin{vmatrix}
\tilde{f}'(\tau)& \tilde{f}''(\tau)&\tilde{f}'''(\tau)\\
g_1'(\tau) & g_1''(\tau)&g_1'''(\tau)\\
g_2'(\tau) & g_2''(\tau)&g_2'''(\tau)\\
\end{vmatrix} d\tau
\end{equation}
to converge for any $t \in \mathbb{R}$ provided~$\gamma$ has negative torsion in a neighborhood of~$-\infty$ (and a similar condition with~$\mathrm{R}$ replaced by~$\mathrm{L}$ and with~$-\infty$ replaced by~$+\infty$ provided~$\gamma$ has positive torsion in a neighborhood of~$\infty$). 

{\bf Claim: under Conditions~\eqref{FirstCondition},~\eqref{ThirdCondition},~\eqref{OneMoreFirstCondition},~\eqref{OneMoreSecondCondition},~\eqref{OneMoreThirdCondition}, and the mentioned convergence conditions for the integrals~\eqref{BulkyIntegral} we can solve Problem~\ref{MainProblem}}. 

As in~\cite{IOSVZ}, by ``solution'' we mean an expression for the function~$\Bell$, which may include roots of implicit equations, differentiations, and integrations. Though at the first sight, the benefit of such a ``solution'' may seem questionable, it occurs to be useful if one has a specific domain~$\Omega$ and a function~$f$ at hand, see examples in the papers~\cite{IOSVZ, ISVZ}, the whole paper~\cite{Vasyunin5} that treats the cases of  functions~$f$ extremely difficult from an algebraic point of view, and other papers on the subject.

It appears that to solve Problem~\ref{MainProblem}, one has to reformulate reasonings from~\cite{IOSVZ} and~\cite{ISVZ} in geometric terms and observe that in such terms they work for a more general setting of the problem considered. For example, the integral~\eqref{BulkyIntegral} plays the role of the force function coming from~$-\infty$ (see~\cite{IOSVZ} for the definition in the case of~$\BMO$) in the general setting. However, the geometric essence of the matter is even more revealed in the example of the chordal domain. We remind the reader that a chordal domain is a type of foliation (see~\cite{IOSVZ} for the definition) that consists of chords, i.e. segments that connect two points of~$\partial \Omega_0$. In the case of the parabolic strip~$g(t) = (t,t^2)$, the chordal domain could match~$\Bell_f$ if and only if it satisfied the cup equation 
\begin{equation*}
\frac{\tilde{f}(b) - \tilde{f}(a)}{b-a} = \frac{\tilde{f}'(b) + \tilde{f}'(a)}{2}; \quad \hbox{$(a,a^2)$ and $(b,b^2)$ are the endpoints of a chord},
\end{equation*} 
and two special differential inequalities (``inequalities for the differentials'') for each of its chord. In the general setting of Problem~\ref{MainProblem}, the cup equation turns into
\begin{equation*}
\begin{vmatrix}
g_1'(a) & g_2'(a) & \tilde{f}'(a)\\
g_1'(b) & g_2'(b) & \tilde{f}'(b)\\
g_1(b)-g_1(a) & g_2(b)-g_2(a) & \tilde{f}(b)-\tilde{f}(a)\\
\end{vmatrix}=0; \quad \hbox{$g(a)$ and $g(b)$ are the endpoints of a chord},
\end{equation*}
which has the following geometrical meaning: the tangent vectors to the curve~$\gamma(t) = (g_1(t),g_2(t),\tilde{f}(t))$ at the points~$a$ and~$b$ lie in one two-dimensional plane with the vector~$\gamma(a) - \gamma(b)$. The special differential inequalities (the so-called inequalities for the differentials) can also be re-stated in purely geometric terms (the triple product of~$\gamma'(a)$,~$\gamma(b) - \gamma(a)$, and the normal to~$\gamma$ at the point~$a$ should be negative; the same should be fulfilled with~$a$ and~$b$ interchanged) and then generalized to fit Problem~\ref{MainProblem}.

In~\cite{IOSVZ} the roots of~$\tilde{f}'''$ played the main role. Indeed, the cups sit on the points where~$\tilde{f}'''$ changes its sign from~$+$ to~$-$. In the general case, the function~$\tilde{f}'''$ should be replaced by the torsion of the curve~$\gamma$. One can see the traces of the torsion in formula~\eqref{BulkyIntegral}. Moreover, now we see that Condition~\eqref{OneMoreSecondCondition} is a straightforward generalization of the regularity condition from~\cite{IOSVZ}.

We recall that in~\cite{IOSVZ} the problem was treated not in the full generality (we assumed that the roots of~$\tilde{f}'''$ were well separated). This narrowed the list of local types of foliations. However, without such an assumption, the collection of figures is wider, see the forthcoming paper~\cite{ISVZ} for the general theory, and the example~\cite{Vasyunin5}, where almost all figures from the general case appear. The latter paper also highlights the notation that becomes very important when there are lots of different figures (it appeared that a foliation corresponds to a special weighted graph). We only mention that all the figures are transferred to the general setting of Problem~\ref{MainProblem}, as well as all the monotonicity lemmas for forces and tails (see~\cite{IOSVZ} for definitions). However, in the general case there are some subtleties concerning different parametrizations of the curves~$g$ and~$\gamma$. To formulate a right analog of a certian monotonicity lemma, one has to choose the right parametrization for it: sometimes it is more convenient to work in the natural parametrization of~$g$, sometimes that of~$\gamma$, sometimes it is useful to lay~$g_1(t) = t$.

Paata Ivanisvili
\smallskip

Department of Mathematics, Michigan State University, East Lansing, MI
48823, USA.
\smallskip

ivanishvili dot paata at gmail dot com
\medskip

Nikolay~N.~Osipov
\smallskip

St. Petersburg Department of Steklov Mathematical Institute RAS,
Fontanka 27, St. Petersburg, Russia;

Norwegian University of Science and Technology (NTNU), IME Faculty,
Dep. of Math. Sci., Alfred Getz' vei 1, Trondheim, Norway.
\smallskip

nicknick at pdmi dot ras dot ru
\medskip

Dmitriy~M.~Stolyrov
\smallskip

Chebyshev Laboratory, St. Petersburg State University, 14th Line, 29b, Saint Petersburg, 199178 Russia;

St. Petersburg Department of Steklov Mathematical Institute RAS,
Fontanka 27, St. Petersburg, Russia.
\smallskip

dms at pdmi dot ras dot ru
\medskip

Vasily~I.~Vasyunin
\smallskip

St. Petersburg Department of Steklov Mathematical Institute RAS,
Fontanka 27, St. Petersburg, Russia.
\smallskip

vasyunin at pdmi dot ras dot ru
\medskip

Pavel~B.~Zatitskiy
\smallskip

Chebyshev Laboratory, St. Petersburg State University, 14th Line, 29b, Saint Petersburg, 199178 Russia;

St. Petersburg Department of Steklov Mathematical Institute RAS,
Fontanka 27, St. Petersburg, Russia.
\smallskip

paxa239 at yandex dot ru
	
\end{document}